\documentclass{article}
\usepackage[utf8]{inputenc}

\usepackage{amsmath}
\usepackage[margin=1in]{geometry}
\usepackage{braket}
\usepackage{stmaryrd}
\usepackage{xcolor}
\usepackage{graphicx}
\usepackage{subfig}
\usepackage{url}

\newcommand{\ymom}{p}
\newcommand{\order}{{\mathcal O}}
\newcommand{\mean}[1]{\braket{#1}}
\newcommand{\fluctint}[1]{\llbracket #1 \rrbracket}
\newcommand\fluct[1]{\left\{{#1}\right\}}
\newcommand{\avg}[1]{\braket{#1}}

\newcommand{\etabar}{\overline{\eta}}
\newcommand{\ubar}{\overline{u}}
\newcommand{\qbar}{\overline{q}}

\usepackage{anyfontsize}

\DeclareMathOperator{\sech}{sech}

\newcommand{\revA}[1]{{#1}}
\newcommand{\revB}[1]{{#1}}
\newcommand{\revAB}[1]{{#1}}
\newcommand{\correct}[1]{{#1}}

\begin{document}


\title{A dispersive effective equation for transverse propagation of planar shallow water waves over periodic bathymetry}
\author{
David I. Ketcheson\thanks{{\texttt{david.ketcheson@kaust.edu.sa}}, Applied Mathematics and Computational Science, CEMSE Division, King Abdullah University of Science and Technology (KAUST), Thuwal, 23955-6900, Kingdom of Saudi Arabia}
\and 
Giovanni Russo\thanks{\texttt{giovanni.russo1@unict.it}, Department of Mathematics and Computer Science, University of Catania, Viale A.\ Doria 6, 95125 Catania, Italy}}

\maketitle

\abstract{We study the behavior of shallow water waves propagating over bathymetry that varies periodically in one direction and is constant in the other.  Plane waves traveling along the constant direction are known to evolve into solitary waves, due to an effective dispersion.  We apply multiple-scale perturbation theory to derive an effective constant-coefficient system of 
equations, \revB{showing that the transversely-averaged wave approximately satisfies a
Boussinesq-type equation, while the lateral variation in the wave is related to certain integral functions of the bathymetry.}
Thus the homogenized equations
not only accurately describe these waves but also predict their full two-dimensional shape in some detail.  Numerical experiments confirm the good agreement between the effective equations and the variable-bathymetry shallow water equations.}

\section{Model Equations and Assumptions}

\begin{figure}
    \center
    \includegraphics[width=0.7\textwidth]{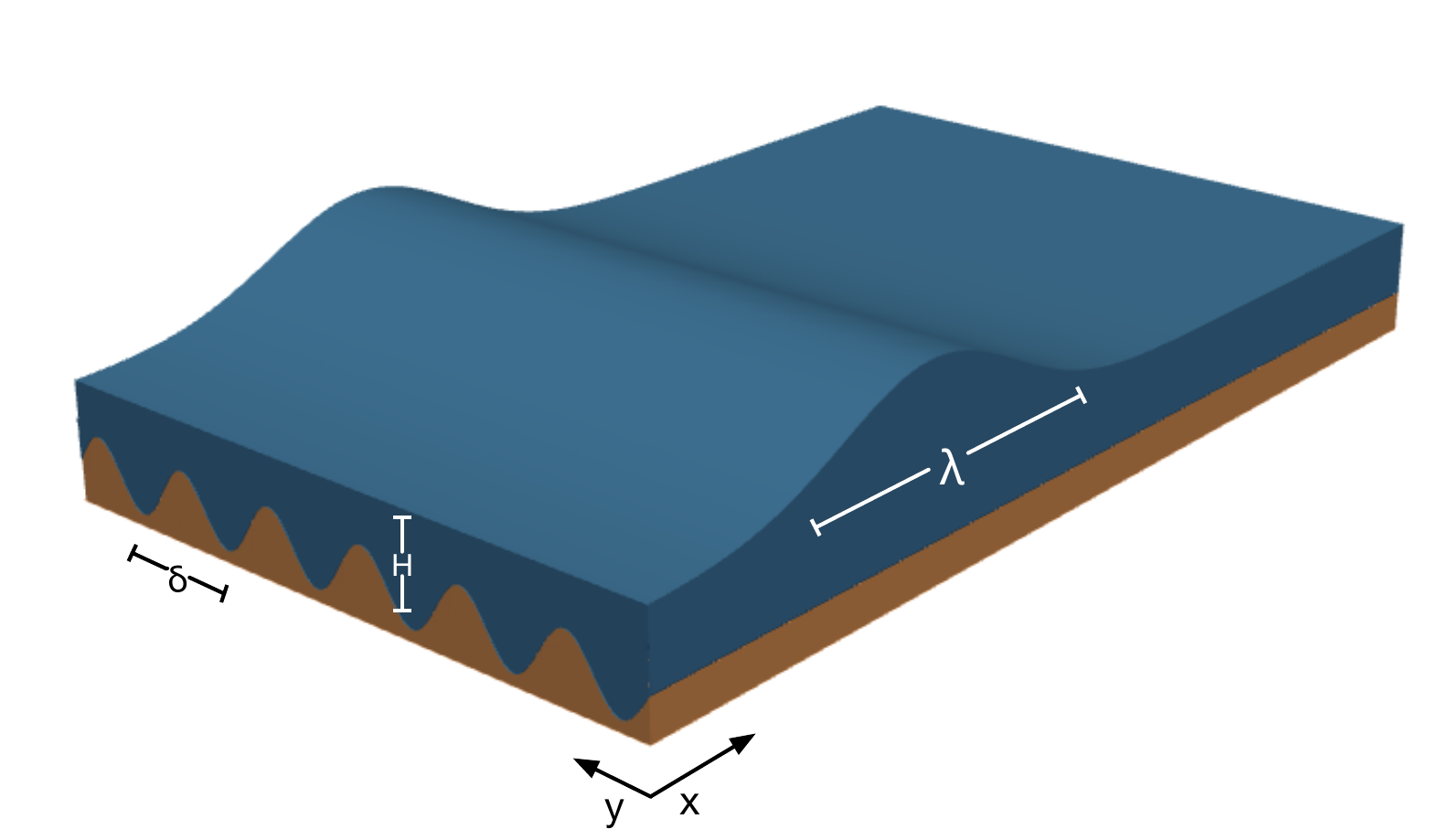}
    \caption{Geometry of the problem studied herein.  The bathymetry $b(y)$ is shown in brown and repeats periodically with period $\delta$ in the $y$ direction.  The regime studied herein is that for which $\lambda \gg \delta$.
    \label{fig:geometry}}
\end{figure}

In this work we study the shallow water wave (or Saint-Venant) model:
\begin{subequations} \label{eq:sw1}
\begin{align}
    h_t + (hu)_x  + (hv)_y & = 0 \\
    (hu)_t + \left( hu^2 + \frac{1}{2}gh^2 \right)_x + (huv)_y & = - g h b_x \\
    (hv)_t + \left( hv^2 + \frac{1}{2} gh^2\right)_y + (huv)_x & = - g h b_y
\end{align}
\end{subequations}
where $g=9.81$ is the gravitational acceleration, $h(x,y,t)$ denotes the depth, $u(x,y,t), v(x,y,t)$ the horizontal 
velocity components, and $b(x,y)$ the bottom elevation (bathymetry).  
As illustrated in Figure \ref{fig:geometry},
we are interested in the behavior of waves propagating over bathymetry that does not depend on $x$, and 
is periodic in $y$ with period $\delta$:
$$
    b(y+\delta) = b(y).
$$
We focus on propagation of initially-planar waves traveling parallel to the $x$-axis:
\begin{align}
    \eta(x,y,0) & = \eta_0(x) & u(x,y,0) & = u_0(x) & v(x,y,0) = 0.
\end{align}
Here $\eta=h+b$ is the surface elevation.
Due to symmetry, this can equivalently be seen as a model for waves in a non-rectangular channel with frictionless walls \cite[Section~1.2]{2021_solitary}, a problem which has also been studied (using other water wave models) for instance in \cite{peregrine1968long,teng1997effects,chassagne2019dispersive}.

It has been shown that linear plane waves traveling parallel to variations in the medium of propagation exhibit effective dispersion if there is variation in the sound speed \cite{quezada_dispersion}.  For nonlinear shallow water waves this effective dispersion 
can lead to the formation of solitary waves even though the equations themselves are
non-dispersive \cite{2021_solitary}.
In the latter work, a partially-heuristic constant-coefficient KdV-type model was shown to approximate the behavior of such waves.  We refer also to \cite{chassagne2019dispersive,gavrilyuk2024geometrical} for development of a similar model and comparison with experiments.
In the present work we develop a more accurate and detailed effective model for these waves.  We apply multiple-scale perturbation analysis to show that these waves 
are described, to leading order, by a Boussinesq-type system with dispersive coefficient depending on the bathymetry.  The effective model describes the full two-dimensional
structure of the waves, and is shown to be in agreement with detailed numerical simulations.
The main novel contributions of this work are:
\begin{itemize}
    \item A constant-coefficient 1D model equation whose solution accurately approximates the average of the 2D variable-coefficient problem;
    \item a computational exploration of solutions of the problem using both finite volume and pseudospectral methods;
    \item analysis and computation of the solitary waves that naturally arise as typical solutions in this problem;
    \item an analytic approximation to those solitary waves.
\end{itemize}

The perturbation approach used here is based on that developed by Yong and coauthors \cite{yong2002solving,leveque2003}.
We have conducted a similar analysis for plane waves propagating perpendicular to the
bathymetric variation; in that case the problem can be reduced to one horizontal dimension 
\cite{ketcheson2023multiscale}.
In the one-dimensional setting, effective dispersion is caused by wave reflection, whereas in the setting of the present work it is caused by propagation \revB{perpendicular to the direction of bathymetric variation and may be described as the result of diffraction \cite{quezada_dispersion}.  It is also reasonable to describe this behavior as the result of refraction, since the deep and shallow regions have different characteristic wave speeds.}
\revB{Compared to the model derived in \cite{2021_solitary}, the model derived herein provides a more detailed (two-dimensional) and accurate description of solutions, in particular for larger waves and longer times.  In principle the technique used here could be carried out to higher orders in order to obtain an even more accurate description.}

Throughout the paper we use dimensional quantities with SI units, so lengths are measured in meters and time in seconds.
The code to reproduce the calculations and figures in this work is available online\footnote{\url{https://github.com/ketch/Shallow_water_2D_homogenization_RR}}.

The rest of the paper is organized as follows.  In Section
\ref{sec:analysis} we perform a multiple-scale analysis
leading to an effective medium equation for the waves of
interest; the main result is equation \eqref{homogenized-system},
which describes the evolution of such waves after averaging
over the $y$-dimension.  In Section \ref{sec:numerical}
we compare solutions of the effective equations with those
of the original variable-bathymetry system \eqref{eq:sw1}.
In Section \eqref{sec:shape} we investigate the shape of
these two-dimensional solitary waves, comparing the predictions
of the multiple-scale analysis with the results of numerical
experiments.
Some conclusions are provided in Section \ref{sec:conclusion}.

\section{Multiple-Scale Analysis} \label{sec:analysis}
The choice of primary variables is a key decision in perturbation analysis of
systems like the one considered here.
One usually works with the conserved variables
$(h, hu, hv)$ in order to include weak solutions,
but here we are interested in strong solutions.
Since we seek to derive a system of equations describing the variation of the solution over long length scales,
we prefer to use quantities that do not necessarily vary on the periodic microscale for near-equilibrium
solutions (see e.g. \cite{leveque2003,ketcheson2023multiscale} for other examples).  In the present setting this indicates that one should use the surface
elevation $\eta$ rather than $h$ and the $y$-momentum $\ymom=hv$ rather than $v$.
After some trial and error we found it best to rewrite \eqref{eq:sw1} in terms of
$(\eta,u,\ymom)$:
\begin{subequations} \label{eq:eta-u-p1}
\begin{align}
    \eta_t + (u(\eta-b))_x + \ymom_y & = 0 \label{eq:2a1} \\
    u_t + u u_x + g \eta_x + \frac{\ymom}{\eta-b} u_y & = 0 \\
    \ymom_t + \left( \frac{\ymom^2}{\eta-b}\right)_y + g(\eta-b)\eta_y + (\ymom u)_x & = 0. \label{eq:2b1}
\end{align}
\end{subequations}
Note that it is not important to write the equations in conservation form since we
are primarily interested in strong solutions of this system.

A second critical choice is that of a small parameter.
We assume that the wavelength of the typical waves we are interested in is long relative to the period $\delta$ of the bathymetry.  
We perform a change of variables, by introducing $\tilde{y} = y/\delta$, so now $b(y) = b(\tilde{y}\delta) = \tilde{b}(\tilde{y})$,
with $\tilde{b}$  a 1-periodic function, and $\partial/\partial y = \delta^{-1}\partial/\partial \tilde{y}$.
We next rewrite the equations in the new coordinates 
$(x,\tilde{y},t)$ and suppress the tildes to obtain
\begin{subequations} \label{eq:eta-u-p}
\begin{align}
    \eta_t + (u(\eta-b))_x + \delta^{-1} \ymom_y & = 0 \label{eq:2a} \\
    u_t + u u_x + g \eta_x + \delta^{-1} \frac{\ymom}{\eta-b} u_y & = 0  \label{eq:2b} \\
    \ymom_t + \delta^{-1} \left( \frac{\ymom^2}{\eta-b}\right)_y + \delta^{-1} g(\eta-b)\eta_y + (\ymom u)_x & = 0. \label{eq:2c}
\end{align}
\end{subequations}
We now look for solutions which are small perturbations, of $O(\delta)$, of the lake at rest given by $(\eta,u,\ymom) = (\eta^0,0,0)$.
\revB{We emphasize that $\eta^0$ denotes the unperturbed surface elevation while $\eta_0(x,y)$ denotes the surface elevation at the initial time.
Furthermore since the problem is invariant under a vertical coordinate shift, we can without loss of generality take $\eta^0=0$, which we will do in the numerical experiments.
}

We assume the quantities $\eta, u, \ymom$ can be written as power series in $\delta$ with the form
\begin{subequations} \label{delta-series}
\begin{align} 
    \eta - \eta^0  & = \delta \eta^1(x,y,t) + \delta^2 \eta^2(x,y,t) + \cdots \\
    u & = \delta u^1(x,y,t) + \delta^2 u^2(x,y,t) + \cdots \\
    \ymom & = \delta \ymom^1(x,y,t) + \delta^2 \ymom^2(x,y,t) + \cdots.
\end{align}
\end{subequations}
Here and throughout this section, superscripts on $\eta, u$, and $\ymom$ denote indices of the asymptotic expansion 
\eqref{delta-series}. When needed,  we shall adopt parentheses to denote exponentiation of these quantities. 
All functions are assumed to be periodic in $y$ with period 1.  In what follows, we use $\braket{\cdot}$ to denote quantities that are averaged with respect to $y$, i.e.\ for any function $f(x,y,t)$ we have
\[
    \braket{f} = \int_0^1 f\, dy.
\]
Notice that if $f$ does not depend on $y$ then $f = \braket{f}$, and that $\braket{\cdot}$ commutes with $x$ and $t$ derivatives, 
i.e.\ $\braket{f}_t = \braket{f_t}$ and $\braket{f}_x = \braket{f_x}$.
Throughout the paper we denote by $$H(y)\equiv \eta^0-b(y)$$ the unperturbed water depth. 

Next, we substitute \eqref{delta-series} into \eqref{eq:eta-u-p} and equate terms at each power of $\delta$.

\subsection{$\order(\delta^0)$}
First we collect all terms proportional to $\delta^0$.
The expansion of \eqref{eq:2b} does not contain any such terms.
From the expansion of \eqref{eq:2a} we obtain 
$
    p^1_y=0
$
while \eqref{eq:2c} gives 
\[
    g H(y)\eta^1_y = 0
\]
From these relations we deduce that $p^1$ and $\eta^1$ do not depend on $y$.

\subsection{$\order(\delta^1)$}
Next, collecting terms proportional to $\delta^1$, we obtain
\begin{subequations} \label{eq:d1}
\begin{align}
    \eta^1_t + H u^1_x + \ymom^2_y & = 0 \label{d1a} \\
    u^1_t + g \eta^1_x + \frac{p^1u^1_y}{H} & = 0 \label{d1b} \\
    \ymom^1_t - \frac{(\ymom^1)^2}{H^2} H + g H \eta^2_y & = 0, \label{d1c}
\end{align}
\end{subequations}
We have concluded already that
$\eta^1$ is independent of $y$, and we assume that $u(x,y,0)$ is independent of $y$, so $u_t$ is initially independent of $y$.  Eq. \eqref{d1b} shows that it will remain so for all time.
Now we average these equations with respect to $y$ -- i.e., we integrate \eqref{eq:d1} with respect to $y$ over one period,
noting that the average of any $y$-derivative is zero.
Since $\eta^1$ is independent of $y$, Eq.~\eqref{d1b} implies that $u^1_t$ is also independent of $y$.
Since we assume that $u(x,y,0)$ is independent of $y$ we deduce further that $u^1$ does not depend on $y$, so we obtain
$$
    \mean{u^1_t} + g \mean{\eta^1_x} = 0.
$$
In order to make the notation more uniform, when considering the evolution of averages in $y$, for any quantity $w$ independent of $y$ we shall indicate it as $\mean{w}$ even if $\mean{w}=w$.

Solving \eqref{d1c} for $\eta^2_y$ and averaging the result gives
$$
    -\frac{\mean{H^{-1}}}{g}
    \avg{\ymom^1_t} = 0.
$$
Since we assume $\ymom(x,y,0)=0$,
this implies $\ymom^1 = 0$.  Returning to \eqref{d1c}, this in turn implies 
$\eta^2_y=0$, so $\eta^2 = \avg{\eta^2}$.
Averaging \eqref{d1a} in $y$ gives
\begin{align}
    \avg{\eta^1_t} + \mean{H} \avg{u^1_x} & = 0.
\end{align}
Taking together these averaged equations, we have the system
\begin{subequations} \label{d1-avg}
\begin{align}
    \avg{\eta^1_t} + \avg{H}\avg{u^1_x} & = 0  \label{d1a-avg}\\
    \avg{u^1_t} + g \avg{\eta^1_x} & = 0, \label{d1b-avg}
\end{align}
\end{subequations}
which is simply the linear wave equation with wave speed $c=\sqrt{g\avg{H}}$.  It is interesting to note that
here the average depth appears in the wave speed, whereas the harmonic average appears for waves propagating perpendicular to the lines of constant bathymetry (see \cite{ketcheson2023multiscale}, and
also \cite{quezada_dispersion}).

We can further manipulate \eqref{d1a} to obtain an expression
for $p^2$, the leading-order term in the y-momentum.
Subtracting \eqref{d1a-avg} from \eqref{d1a} 
we obtain 
\[
    -p^2_y = \fluct{\eta^1_t} + \fluct{H}u^1_x 
\]
where, for any function $f$, we denote by 
\[
    \fluct{f} \equiv f - \avg{f}
\]
the {\em fluctuating part\/} of $f$.
Considering that $\eta^1=\avg{\eta^1}$ we have
\begin{equation}\label{eq:p2_0}
    p^2_y = -\fluct{H}u^1_x = 
    \correct{-\fluct{H}\mean{u^1_x}  }
\end{equation}
Integrating Eq.~\eqref{eq:p2_0} and imposing that the two sides have the same average, one gets
\begin{equation} 
p^2(x,y,t) = -\fluctint{H}\correct{\mean{u^1_x}} + \avg{p^2} \label{eq:p2}
\end{equation}
where, for any function of $f(y)$, $\fluctint{f}$ denotes the integral of the fluctuating part of $f$:
$$
    \fluctint{f} = \int_s^y (f(\xi)-\mean{f}) d\xi \ \ \text{where $s$ is chosen so that } \mean{\fluctint{f}} = 0.
$$


\subsection{$\order(\delta^2)$}
Collecting terms proportional to $\delta^2$, we obtain
\begin{subequations} \label{eq:d2}
\begin{align}
    -\ymom^3_y & = \eta^2_t + Hu^2_x + (\avg{\eta^1} \avg{u^1})_x \label{d2a} \\
    0 & = u^2_t + \avg{u^1} \avg{u^1_x} + g \avg{\eta^2_x}  \label{d2b} \\
    -g \eta^3_y & = \frac{1}{H}(\ymom^2_t + g \avg{\eta^1} \avg{\eta^2_y}) = \frac{1}{H}(-\fluctint{H}\avg{u^1_{xt}} + \avg{\ymom^2}_t). \label{d2c}
\end{align}
\end{subequations}
In the last line we have used that $\eta^2$ is
independent of $y$ and Equation \eqref{eq:p2}.
The only term in \eqref{d2b} that could depend on $y$ is $u^2_t$, so it must be independent of $y$; i.e. $u^2 = \avg{u^2}(x,t)$.  Thus we have
\begin{align} \label{d2b-avg}
    \avg{u^2_t} + \avg{u^1} \avg{u^1_x} + g \avg{\eta^2_x} & = 0.
\end{align}
Taking the average of \eqref{d2a} gives
\begin{align} \label{d2a-avg}
    \avg{\eta^2_t} + \mean{H}\avg{u^2_x} + (\avg{\eta^1} \avg{u^1})_x = 0.
\end{align}
Subtracting this from \eqref{d2a} and integrating in $y$, we get
\begin{align} \label{eq:p3}
    \ymom^3(x,y,t) & = -\fluctint{H}\avg{u^2_x} + \avg{\ymom^3}.
\end{align}
For simplicity we now specialize our analysis to bathymetry profiles for which $\mean{H^{-1}\fluctint{H}}=0$,
which holds for instance for the piecewise-constant
or sinusoidal bathymetries studied below (see \cite{ketcheson2023multiscale}, Proposition 5).  Then
taking the average of \eqref{d2c}, we find that $\avg{\ymom^2}_t=0$.
Since $p^2(x,y,0) = 0$, it follows that $\avg{p^2} = 0$.
Then integrating \eqref{d2c} in $y$ yields
\begin{align} \label{eq:eta3}
    \eta^3(x,y,t) & = \frac{1}{g}\fluctint{H^{-1}\fluctint{H}}\avg{u^1_{xt}} + \avg{\eta^3}.
\end{align}

Based on what we have determined
up to this point, we can write the series \eqref{delta-series}
more simply as
\begin{subequations} \label{delta-series-simp}
\begin{align} 
    \eta - \eta^0 & = \delta \eta^1(x,t) + \delta^2 \eta^2(x,t) + \cdots \\
    u & = \delta u^1(x,t) + \delta^2 u^2(x,t) + \cdots \\
    \ymom & = \delta^2 \ymom^2(x,y,t) + \cdots.
\end{align}
\end{subequations}

Let $\overline{\eta} = \delta^{-1}\avg{\eta-\eta^0}$ and $\overline{u}=\delta^{-1}\avg{u}$.
By adding $\delta$ times \eqref{d1a-avg} to $\delta^2$ times \eqref{d2a-avg}, we get an 
approximate equation for the evolution of $\etabar$:
\begin{subequations}
\begin{align}
    \delta\left(\etabar_t + \avg{H}\ubar_x\right) + \delta^2(\etabar \ \ubar)_x & = \order(\delta^3).
\end{align}
Similarly, by adding $\delta$ times \eqref{d1b-avg} to $\delta^2$ times \eqref{d2b-avg}, we get an
approximate equation for the evolution of $\ubar$:
\begin{align}
    \delta\left(\ubar_t + g\etabar_x\right) + \delta^2(\ubar \ \ubar_x) & = \order(\delta^3).
\end{align}
\end{subequations}
We see that up to this order, the y-averaged variables satisfy
a nonlinear first-order hyperbolic system.
We proceed with the analysis at the next order,
where we expect to see dispersive terms.

\subsection{$\order(\delta^3)$}
Collecting terms proportional to $\delta^3$, we obtain
\begin{subequations}
\begin{align}
    \eta^3_t + H u^3_x + (\eta^1u^2 + \eta^2 u^1)_x & = -\ymom^4_y \label{d3a} \\
    u^3_t + (u^1 u^2)_x + \fluctint{H^{-1}\fluctint{H}} u^1_{xtx} + g \eta^3_x & = 0 \label{d3b} \\
    \frac{1}{H} \left( \ymom^3_t - H^{-2}(\ymom^2)^2 H' + 2 H^{-1} \ymom^2 \ymom^2_y + g \eta^1 \eta^3_y + (\ymom^2 u^1)_x \right) & = -g \eta^4_y.  \label{d3c}
\end{align}
\end{subequations}
Averaging \eqref{d3a} yields
\begin{align} \label{d3a-avg}
    \avg{\eta^3}_t + \mean{H u^3}_x + (\avg{\eta^1} \avg{u^2} + \avg{\eta^2} \avg{u^1})_x & = 0.
\end{align}
Since $u^3$ may depend on $y$, we must work directly with the average $\mean{Hu^3}$.  We therefore multiply \eqref{d3b} by $H$ before averaging, to obtain
\begin{align}
    \mean{Hu^3}_t + \mean{H}(\avg{u^1} \avg{u^2})_x - \mu \avg{u^1_{xxt}} + g \mean{H}\avg{\eta^3_x} & = 0.
\end{align}
where
\begin{align}
    \mu & = -\mean{H \fluctint{H^{-1}\fluctint{H}} } = \mean{H^{-1}(\fluctint{H})^2}.
\end{align}
The last equality comes from the general property $\mean{a\fluctint{b}}=-\mean{\fluctint{a}b}$ for all functions $a(y)$ $b(y)$ \cite[Appendix A]{yong2002solving} (here we take $a=H, b=H^{-1}\fluctint{H}$).

Introducing $\avg{q^j} = \mean{Hu^j}$, this is
\begin{align}
    \avg{q^3}_t + \mean{H}^{-1}(\avg{q^1} \avg{q^2})_x - \mean{H}^{-1}\mu \avg{q^1_{xxt}} + g \mean{H}\avg{\eta^3_x} & = 0. \label{d3b-avg}
\end{align}
Here we made use of the fact that $u^1$ and $u^2$ are independent of y.
Averaging \eqref{d3c}, after a number of tedious calculations, yields $\avg{\ymom^3} = 0$.

\subsection{Governing equations for averaged variables}
Let $\qbar = \sum_j \delta^{j-1} \avg{q^j}$.
We now add $\delta\mean{H}$ times \eqref{d1b-avg}, plus $\delta^2\mean{H}$ times \eqref{d2b-avg}, plus
$\delta^3$ times \eqref{d3b-avg}.  This gives
\begin{subequations} \label{homogenized-system}
\begin{align}
    \delta( \qbar_t + g\mean{H}\etabar_x) + \delta^2 \mean{H}^{-1} \qbar \ \qbar_x & = \delta^3 \mean{H}^{-1} \mu \qbar_{xxt} + \order(\delta^4)
\end{align}
Similarly, adding $\delta$ times \eqref{d1a-avg} with $\delta^2$ times \eqref{d2a-avg} with $\delta^3$ times \eqref{d3a-avg} results in
\begin{align}
    \delta(\etabar_t + \qbar_x) + \delta^2 \mean{H}^{-1}(\etabar \ \qbar)_x & = \order(\delta^4).
\end{align}
\end{subequations}
It turns out that this system is identical to the so-called \emph{classical Boussinesq system} that was originally derived as a model for long-wavelength waves over a flat-bottom channel.  Remarkably, here it has arisen in a completely different way, starting from the non-dispersive Saint-Venant system, and with dispersion arising purely from the effect of a non-flat bottom.
In the present context the coefficients of the convective and dispersive terms depend on the bathymetry $b(y)$ and so their relative magnitude can be quite different based on the chosen geometry.
This system is known to be well-posed \cite{schonbek1981existence,amick1984regularity,bona2004boussinesq}.

In principle the asymptotic analysis can be carried out to higher order, deriving additional high-order effective dispersive terms.  Of course, it should be kept in mind that by starting from the Saint-Venant system \eqref{eq:sw1} we have already discarded certain higher-order effects that might compete with or dominate the additional terms obtained through such an analysis.  This will depend on the relative size of the shallowness parameter $H/\lambda$ and the bathymetry parameter $\delta/\lambda$.

\revAB{
\subsection{Two-dimensional wave structure}\label{sec:2d-structure}
In addition to providing an effective model for wave propagation in terms of $y$-averaged quantities, the homogenization process also allows us to determine the variation of solutions with respect to the $y$ coordinate.

For the $y$-momentum, from equations \eqref{eq:p2} and \eqref{eq:p3} together with the fact that $\mean{p}=0$ we immediately obtain that the variation in $y$ is proportional to $\fluctint{H}$:
\begin{align} \label{p-full}
    p(x,y,t) & \approx -\fluctint{H} \avg{u_x}.
\end{align}

From \eqref{d1b-avg} we have that $\mean{u^1_t} \approx -g\mean{\eta^1_x}$,
and therefore $\mean{u^1_{xt}} \approx -g\mean{\eta^1_{xx}}$, so then from \eqref{eq:eta3} we obtain
\begin{align}
    \eta(x,y,t) - \eta^0 & = \delta \mean{\eta^1}(x,t) + \delta^2 \mean{\eta^2}(x,t) + \delta^3(\mean{\eta^3}(x,t) - \fluctint{H^{-1}\fluctint{H}} \mean{\eta^1_{xx}}) + \order(\delta^4),
\end{align}
or equivalently
\begin{align} \label{eta-y-variation}
    \eta(x,y,t) \approx \mean{\eta}(x,t) - \delta^2\fluctint{H^{-1}\fluctint{H}} \mean{\eta_{xx}},
\end{align}
showing that the leading variation with respect to $y$ is proportional to 
$\fluctint{H^{-1}\fluctint{H}} \mean{\eta_{xx}}$.

Finally, from the leading-order {\em linear} system \eqref{d1-avg}, looking for simple waves, 
we have that
\begin{align}
    u(x,y,t) & \approx \pm \sqrt{g/\avg{H}} (\eta - \eta^0) \label{u-full}
\end{align}
with the plus sign for right-going waves
and minus for left-going waves.
}

\section{Numerical comparison}\label{sec:numerical}
In this section we explore the accuracy of the  homogenized approximation by comparing 
its numerical solutions to numerical solutions of the original system \eqref{eq:sw1}. 
We start by discussing the methods adopted for the numerical solution of both the original system 
\eqref{eq:sw1} and the homogenized  system \eqref{homogenized-system}.  

\subsection{Numerical discretization of the homogenized equations}
We solve the homogenized equations \eqref{homogenized-system} with a Fourier
pseudospectral discretization in space and 
explicit 3-stage 3rd-order SSP Runge-Kutta integration in time.
We can write this system as
\begin{subequations} \label{semi-discrete}
\begin{align}
    \etabar_t & = - \qbar_x - \delta \avg{H}^{-1} (\etabar \ \qbar)_x \\
    \qbar_t & = -(1-\delta^2 \avg{H}^{-1} \mu \partial_x^2)^{-1} \left(g \avg{H}\etabar_x + \delta \avg{H}^{-1} \qbar \ \qbar_x \right)
\end{align}
\end{subequations}
We discretize in the standard pseudospectral way and then apply
the inverse elliptic operator $(1-\delta^2 \avg{H}^{-1} \mu \partial_x^2)^{-1}$
in Fourier space, which does not require the solution of any
algebraic system.  
We can therefore integrate the pseudospectral 
semi-discretization of \eqref{semi-discrete} efficiently with an 
explicit Runge--Kutta method.

For the spatial domain, we take $x\in[-L,L]$ where $L$ is chosen
large enough that the waves do not reach the boundaries before
the final time.

\subsection{Numerical methods for the variable-bathymetry shallow water system}
For the solution of the first-order variable-coefficient
hyperbolic shallow water system \eqref{eq:sw1} we use
two different approaches, depending on the nature of the bathymetry.  
Accurate solution of this system is much more expensive as it requires a much finer spatial
mesh, in order to resolve the bathymetric variation and its effects, and it requires the
solution of a problem in two space dimensions.
applied at $x=0$.

For piecewise-constant (discontinuous) bathymetry, we use the finite volume code Clawpack \cite{LeVeque-FVMHP,mandli2016clawpack}, employing the 
\emph{SharpClaw} algorithm, based on 5th-order WENO
reconstruction in space and 4th-order Runge--Kutta integration in time \cite{KetParLev13}.
This algorithm is well adapted to handle the lack
of regularity in both the coefficients and the solution.
For continuous bathymetry, we again use the Clawpack code and we also compute the solution with a standard Fourier collocation pseudospectral method in space and 4th-order Runge--Kutta integration in time.
Ordinarily one would avoid the use of spectral methods for a first-order hyperbolic problem, but since we focus on scenarios in which shocks do not form, this method performs well and is more efficient than a finite volume discretization,
as long as the bathymetry is continuous.

\subsection{Smooth bathymetry}
First we consider the smoothly-varying bathymetry
\begin{subequations} \label{sinusoidal-setup}
\begin{align}
    b & = -1 + \frac{3}{10} \sin(2\pi y) \label{sinbath} \\
    \eta(x,y,0) & = \frac{1}{20} \exp(-(x/5)^2) \\
    u(x,y,0) & = 0 \\
    v(x,y,0) & = 0
\end{align}
\end{subequations}
\revB{Although the same problem is solved using the PS and FV methods, the computational setup is slightly different.

For the PS code we take $(x,y)\in [-1000,1000]\times[-1/2,1/2]$ and impose periodic boundary conditions at all boundaries.  The final time is chosen such that the waves do not reach $x=\pm1000$.

For the FV code we can save computational effort by taking only
the right half of the domain: $(x,y)\in [0,1000]\times[-1/2,1/2]$.  Due to the symmetry of the solution, we can impose a reflecting boundary condition at $x=0$ and obtain the solution of the same problem, restricted to the right half of the domain.
However, this is still quite expensive, so to save even
more computational effort we do as follows.
We take a much smaller domain ($(x,y)\in [0,100]\times[-1/2,1/2]$).  We impose a reflecting boundary condition at $x=0$ initially; once the waves have moved away from the origin  we impose periodic boundary conditions in $x$.  This allows us to simulate the right-going wave train with high resolution at a reasonable computational cost.  The simulation ends long before the leading wave would begin to catch up to the tail of the wave train.
This is possible because the solution is nearly constant far away from the main waves.}

\begin{figure}
    \center
    \includegraphics[width=\textwidth]{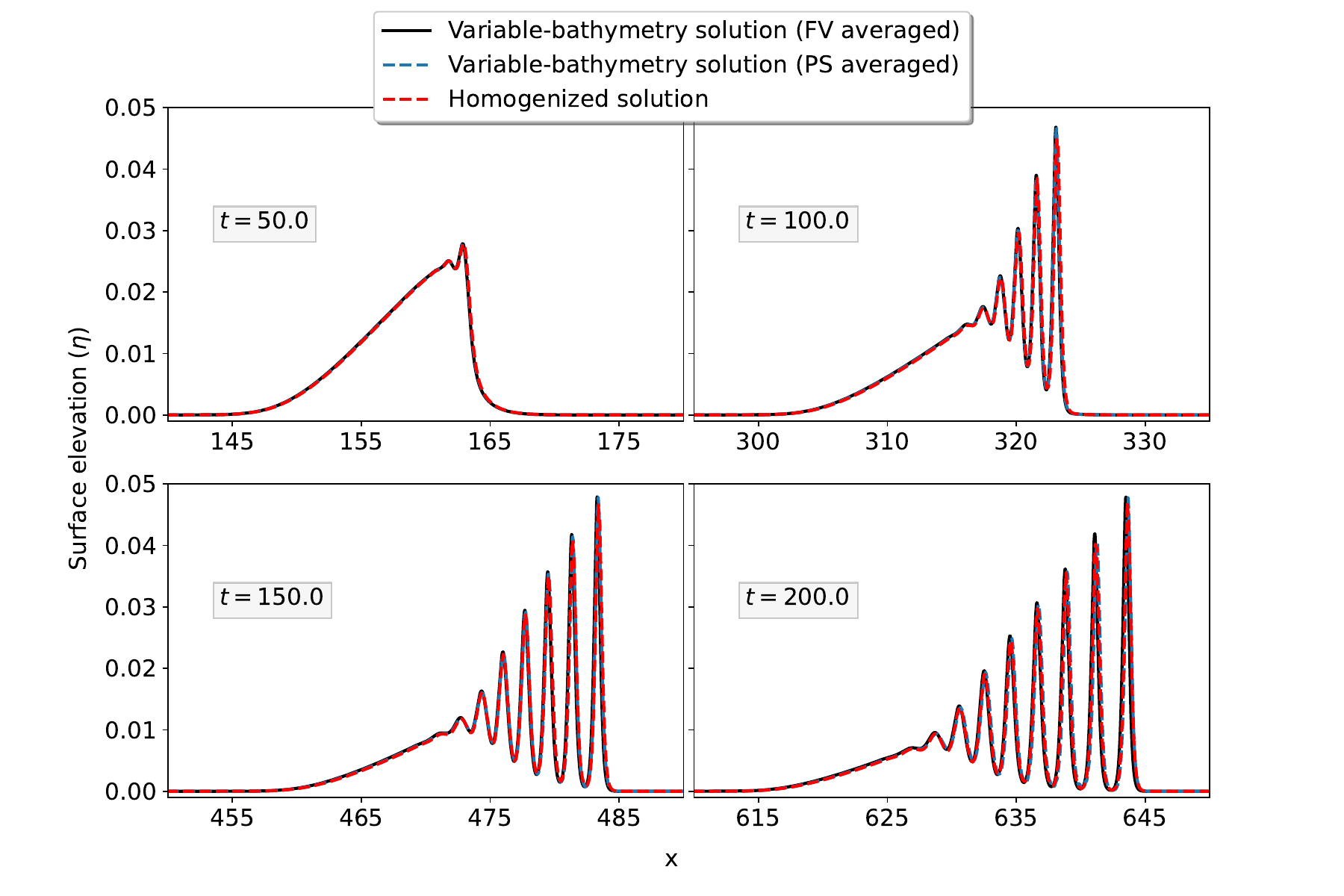}
    \caption{Comparison of homogenized and direct solutions, for sinusoidal bathymetry \eqref{sinbath}.
        The surface elevation $\eta - \eta^0$ is shown.  \label{fig:comparison_sin}}
\end{figure}

  Results are shown in Figure \ref{fig:comparison_sin}. 
  The initial surface perturbation splits into a left-going and right-going pulse, each of which
  eventually resolves into a series of traveling waves.  
  We see a remarkably close agreement between all three solutions, up to $t=200$.

  To obtain the results shown here, we used a mesh of $32000 \times 32$ points for the 2D PS code and $16000\times 160$ points for the 2D FV code.  For the 1D homogenized equations, the pseudospectral simulation was performed on a grid with 32000 points.

\subsection{Piecewise-constant bathymetry}
We next consider the discontinuous bathymetry:

\begin{align} \label{pwc-setup}
    b(x,y) & = \begin{cases} -2/5 & -1/2 \le y < 0 \\ -8/5 & 0 \le y < 1/2 \end{cases}
\end{align}
with the same initial data as in \eqref{sinusoidal-setup}.
In this case we cannot use the 2D PS solver due to the lack of continuity of the solution.
For the FV simulation, the domain and boundary conditions are set up in the same way as for the problem above.

\begin{figure}
    \center
    \includegraphics[width=\textwidth]{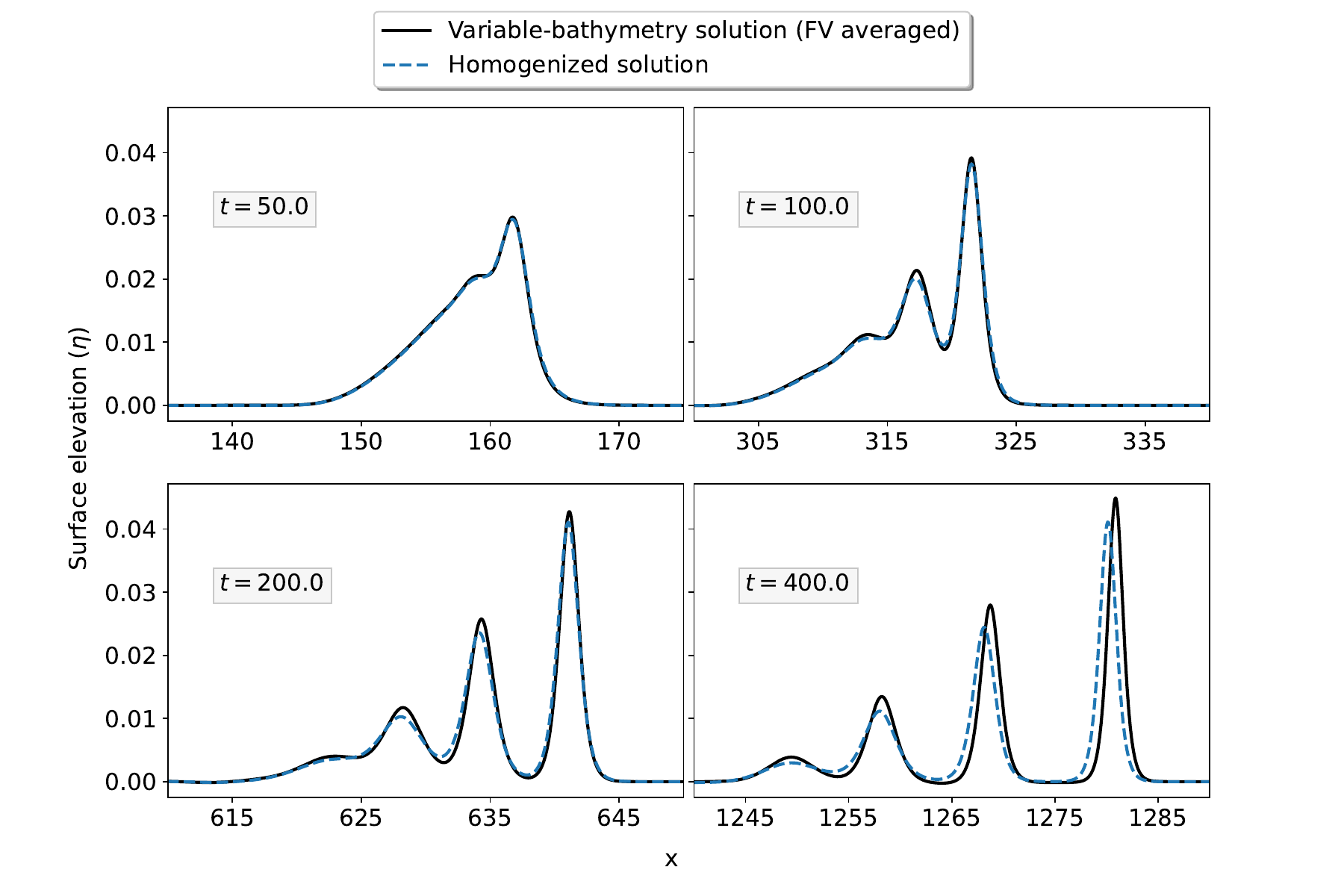}
    \caption{Comparison of homogenized and direct solutions, for piecewise-constant bathymetry \eqref{pwc-setup}.
        The surface elevation $\eta - \eta^0$ is shown.  \label{fig:comparison}}
\end{figure}

In Figure \ref{fig:comparison} we
show snapshots of the right-going pulse.  We see extremely close agreement between the
solutions, with some differences visible at late times, after the pulse has propagated
for hundreds of meters. 

\begin{figure}
\center
    \includegraphics[width=0.75\textwidth]{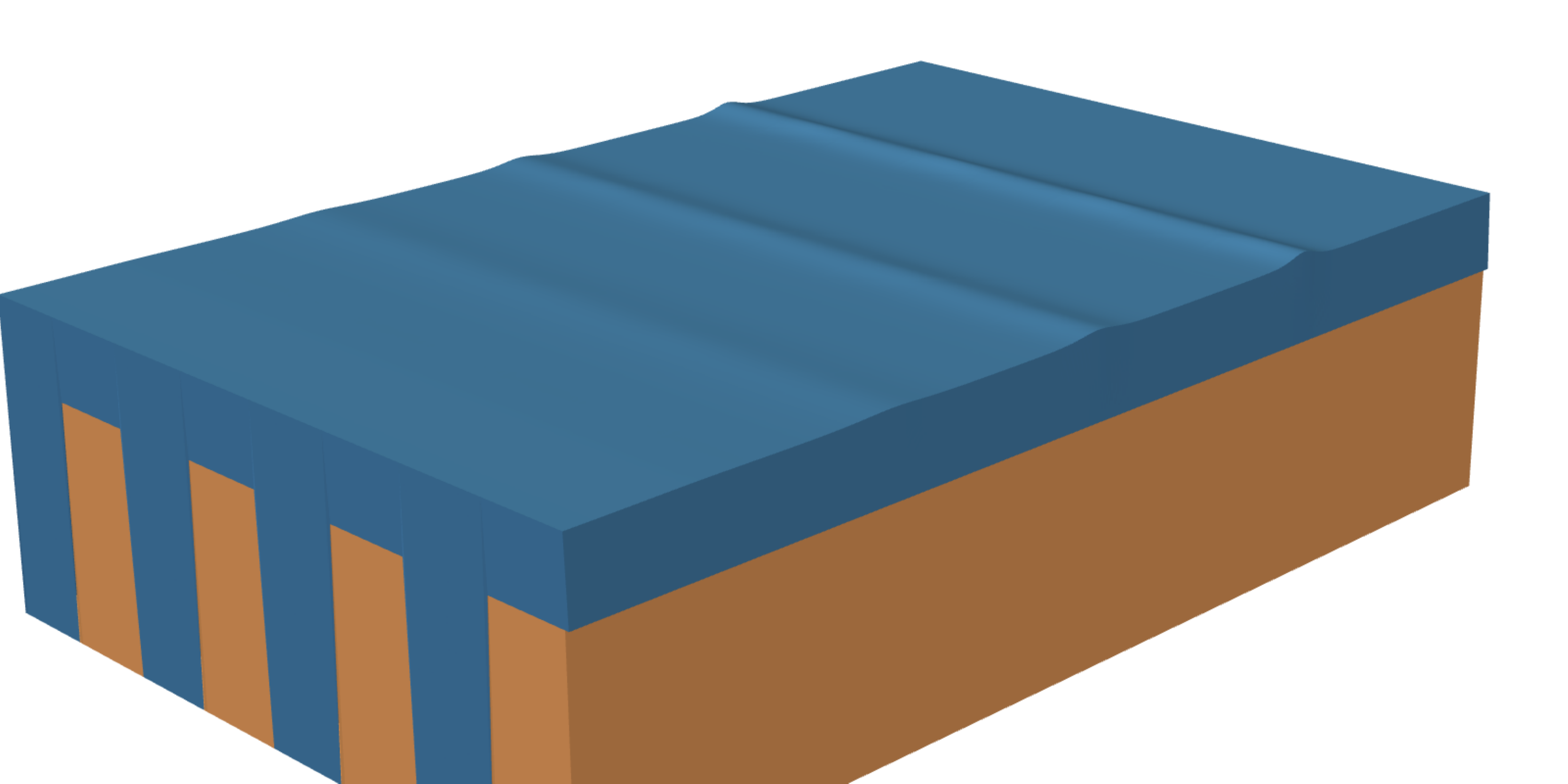}
    \caption{\revB{Three-dimensional rendering of the FV solution shown in the last
    panel of Figure \ref{fig:comparison}.  Results shown are to scale, except that the $x$-axis has been compressed by a factor of 8 to make it easier to see the waves.}}
    \label{fig:iso2}
\end{figure}

\begin{figure}
\center
    \includegraphics[width=0.75\textwidth]{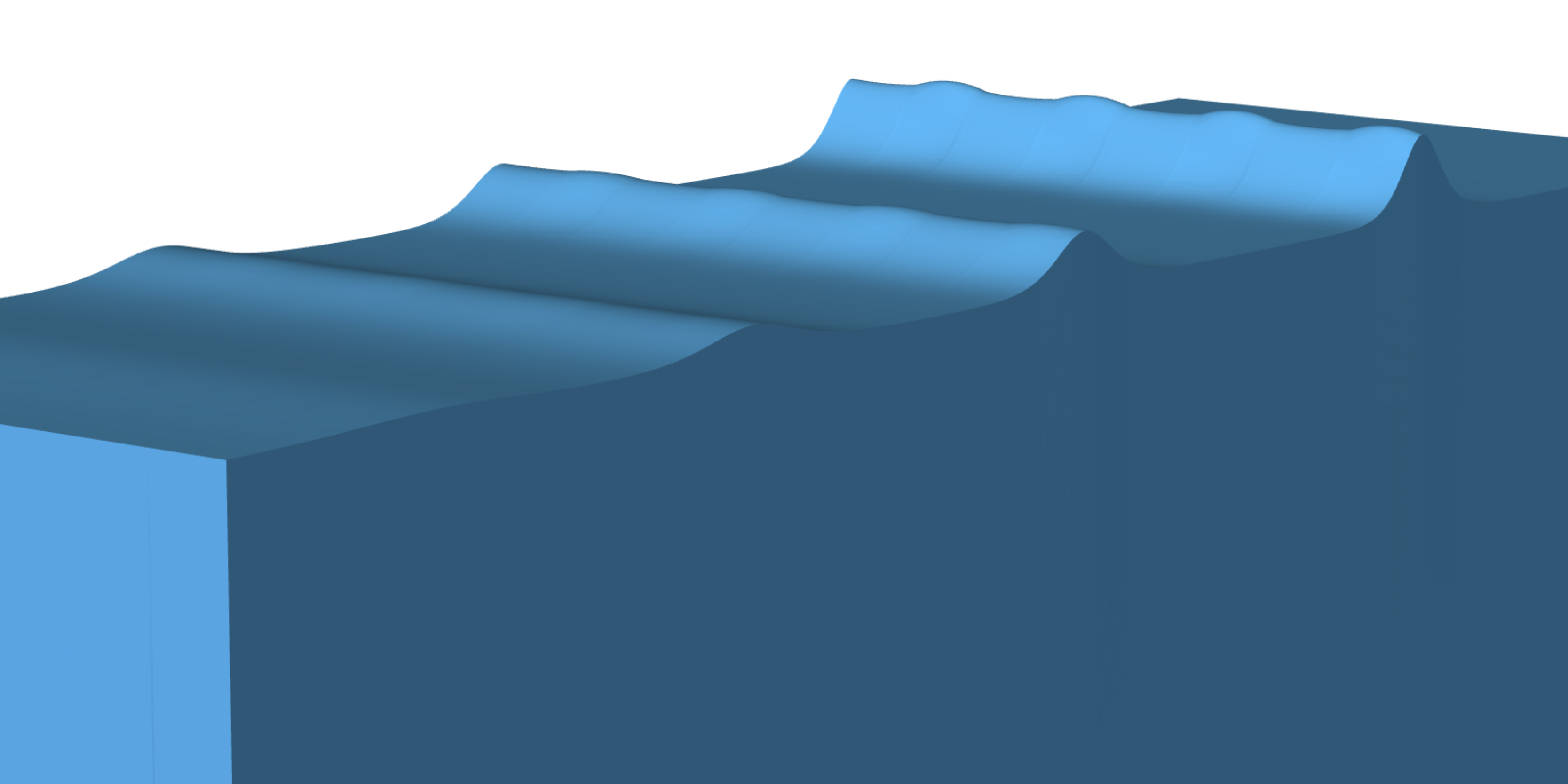}
    \caption{\revB{Closeup of the surface waves shown in Figure \ref{fig:iso2}.
    The vertical variation has been exaggerated by 10x and the $x$-axis has been compressed by a factor of 8 to make it easier to see the waves.
    }}
    \label{fig:iso_zoom}
\end{figure}

\begin{figure}
\center
    \includegraphics[width=0.75\textwidth]{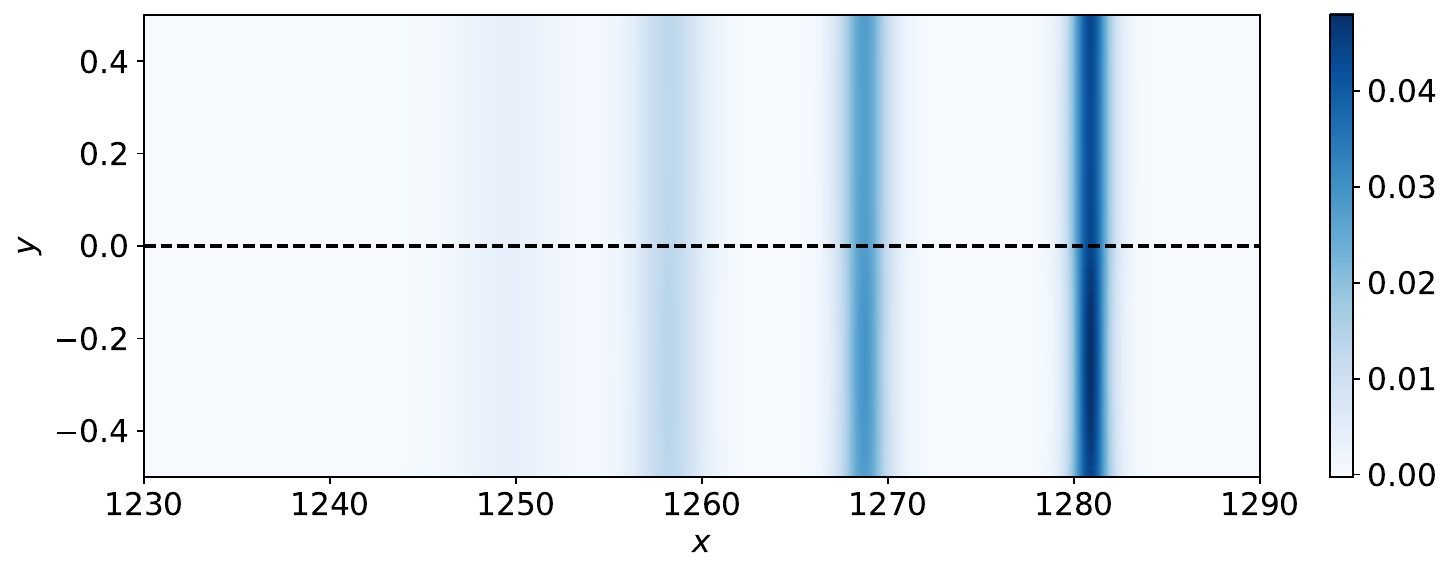}
    \caption{\revB{Surface elevation for the FV solution shown in the last panel of Figure \ref{fig:comparison}.}}
    \label{fig:surface_pcolor}
\end{figure}

\begin{figure}
\center
    \includegraphics[width=0.75\textwidth]{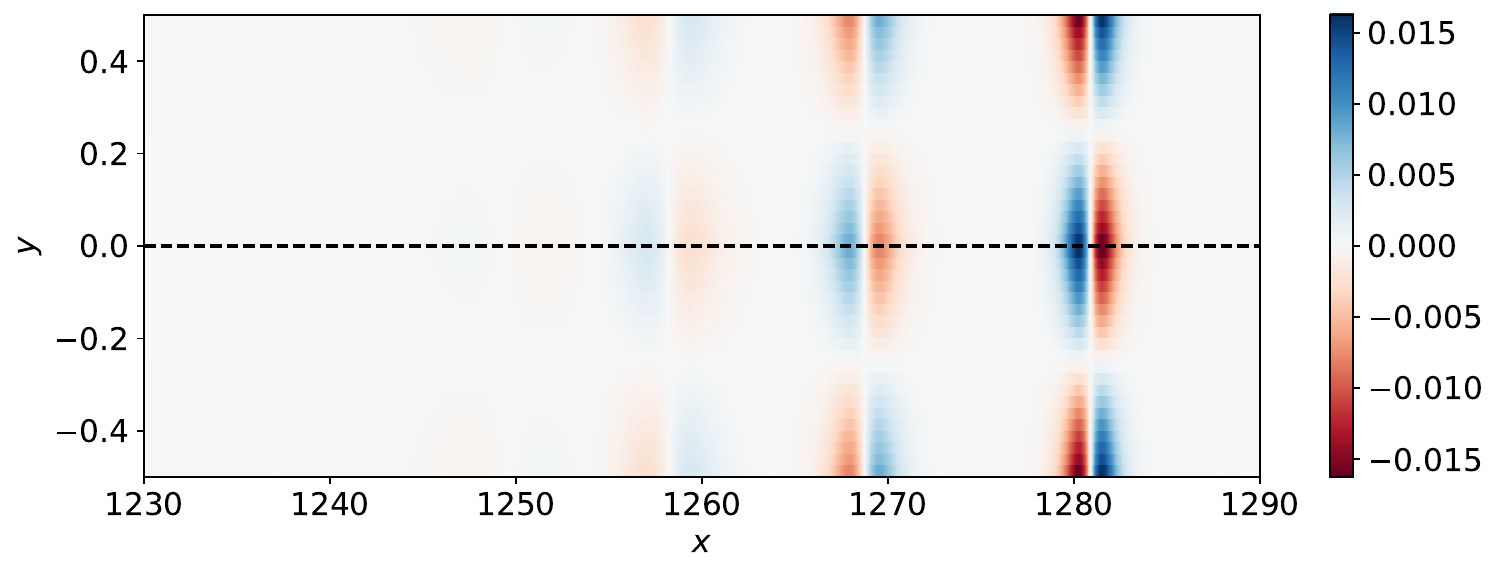}
    \caption{\revB{$y$-momentum for the FV solution shown in the last panel of Figure \ref{fig:comparison}.}}
    \label{fig:hv_pcolor}
\end{figure}

\revB{
The solution at the final time is also shown in Figures \ref{fig:iso2}-\ref{fig:hv_pcolor}.  Although the characteristic speed varies greatly as a function of $y$, the surface elevation profiles remain almost perfectly planar.  The mechanism for this can be seen in Figure \ref{fig:hv_pcolor}, which exhibits a small flow from the deep region to the shallow region at the front of
each solitary wave, and a similar flow from the shallow region to the deep region at the back of each solitary wave.
One can also view each solitary wave as a superposition of an upward-traveling wave with an oscillatory shape (shown in blue in Figure \ref{fig:hv_pcolor}) and a downward-traveling wave with an oscillatory shape (shown in red).
These two waves combine to yield highly planar surface and $x$-velocity fields.

As predicted by equation \eqref{eta-y-variation}, the
wave height does vary to a small degree with $y$;
this can be seen in Figure \ref{fig:iso_zoom}.
We examine this variation
more carefully in Section \ref{sec:shape}.}

\section{Solitary wave shape} \label{sec:shape}
In this section we study the shape of the solitary waves observed in numerical simulations and compare them with predictions based on the homogenized equations.  We first consider traveling wave solutions of the homogenized equations, and then investigate the full two-dimensional solitary waves in more detail.

\subsection{Traveling wave solutions of the homogenized equations}
Now we consider the problem of finding the traveling wave solution for the homogenized system \eqref{homogenized-system}.
Neglecting the higher order term, after dividing by $\delta$, and neglecting terms of $\order{(\delta^3)}$, the system can be written in the form
\begin{align}
    q_t + a_1 \eta_x + a_2 q q_x - \tilde{\mu}q_{xxt} & = 0\\
    \eta_t + q_x + a_2(\eta q)_x & = 0
\end{align}
where we set $a_1:=g\mean{H}$, $a_2 := \delta/\mean{H}$, $\tilde{\mu} := \delta^2\mu/\mean{H}$.

\revA{This is the classical Boussinesq model with linear dispersion, which has been widely studied in the literature (see, for example, \cite{chen2000solitary}). To be more self-contained we summarize here the main steps of the procedure.}

We look for traveling waves which depend only on $\xi = x-Vt$, propagating on a lake at rest, so that the unperturbed state is $q_0=0$, and, with a suitable choice of the frame of reference, $\eta_0 = 0$. Here $V$ is the traveling speed of the wave. Assuming $\eta$ and $q$ are functions of $\xi$, we obtain that the wave \revA{satisfies an equation of the form}
\begin{equation}
    \label{eq:q''}
    q'' = F(q),
\end{equation}
with
\[
    \eta = \frac{q}{V-a_2 q}.
\]
Multiplying by $q'$ and integrating we obtain 
\revA{the analogue of total energy conservation}
\begin{align} \label{qp2}
    \frac12 (q')^2 + U(q) = E,
\end{align}
with 
\[
    U(q) = \left(\frac16a_2q^3-\frac12Vq^2-\frac{a_1}{a_2}q
    -\frac{a_1}{a_2^2}V\log(1-a_2 q/V)\right)/(\tilde{\mu} V).
\]
The trajectories of the material point in phase space $(q,\dot{q})$ are the lines which maintain constant total energy. 
Notice that if we approximate the $\log$ term by the first term in its expansion about $q=0$, then the solution of \eqref{qp2} is a hyperbolic secant squared.  Thus we expect that solitary waves will be close to this shape.

\begin{figure}
    \centering
    \includegraphics[width=0.3\linewidth]{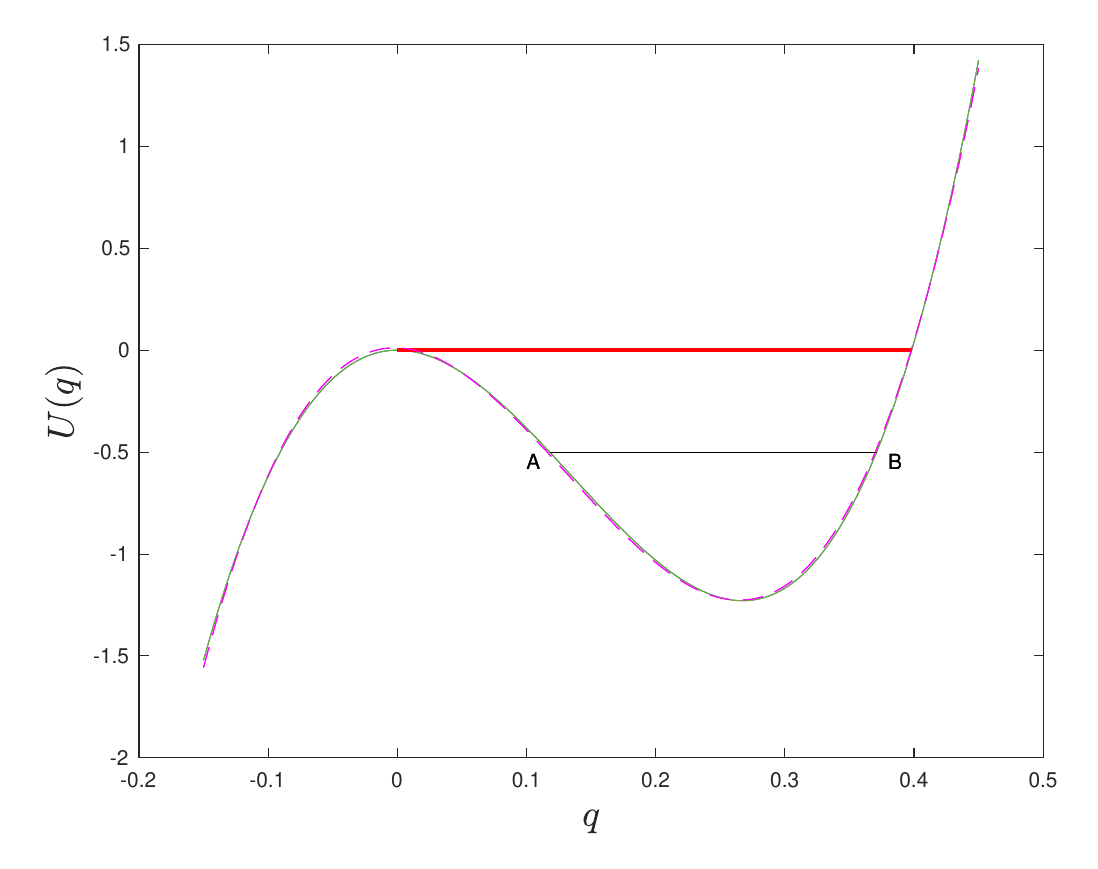}
    \includegraphics[width=0.3\linewidth]{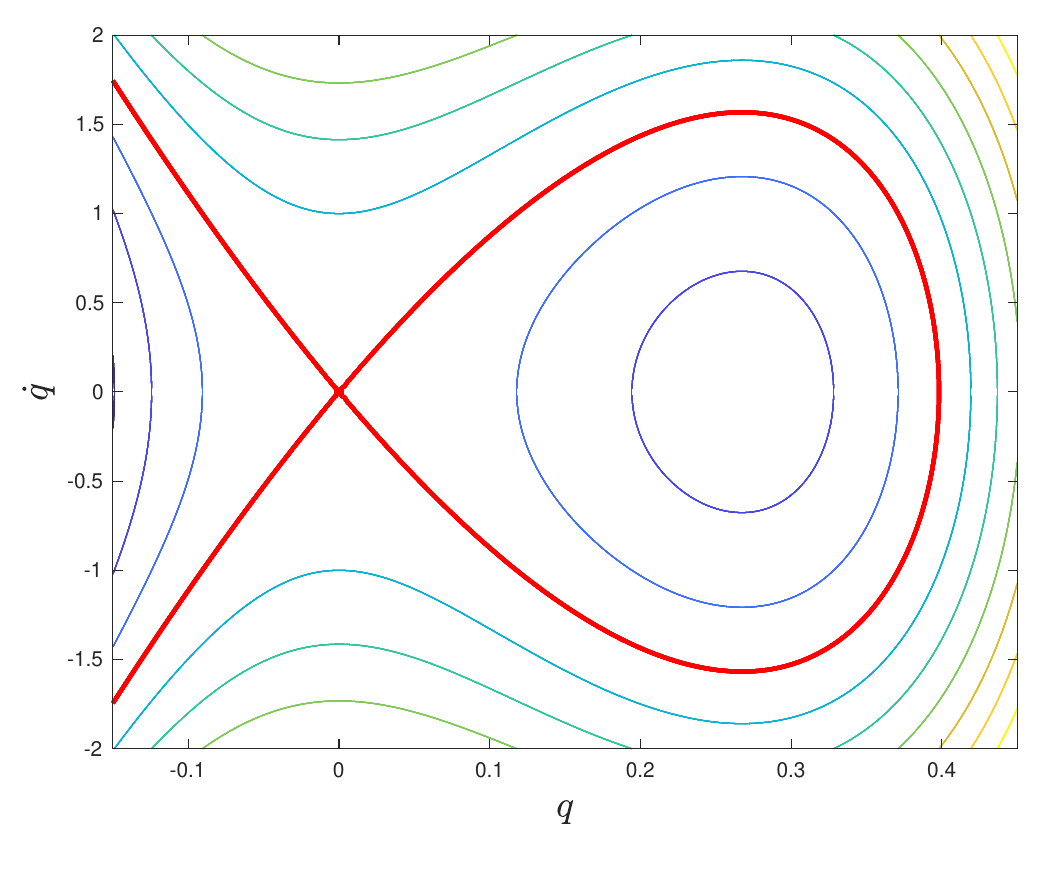}
    \includegraphics[width=0.3\linewidth]{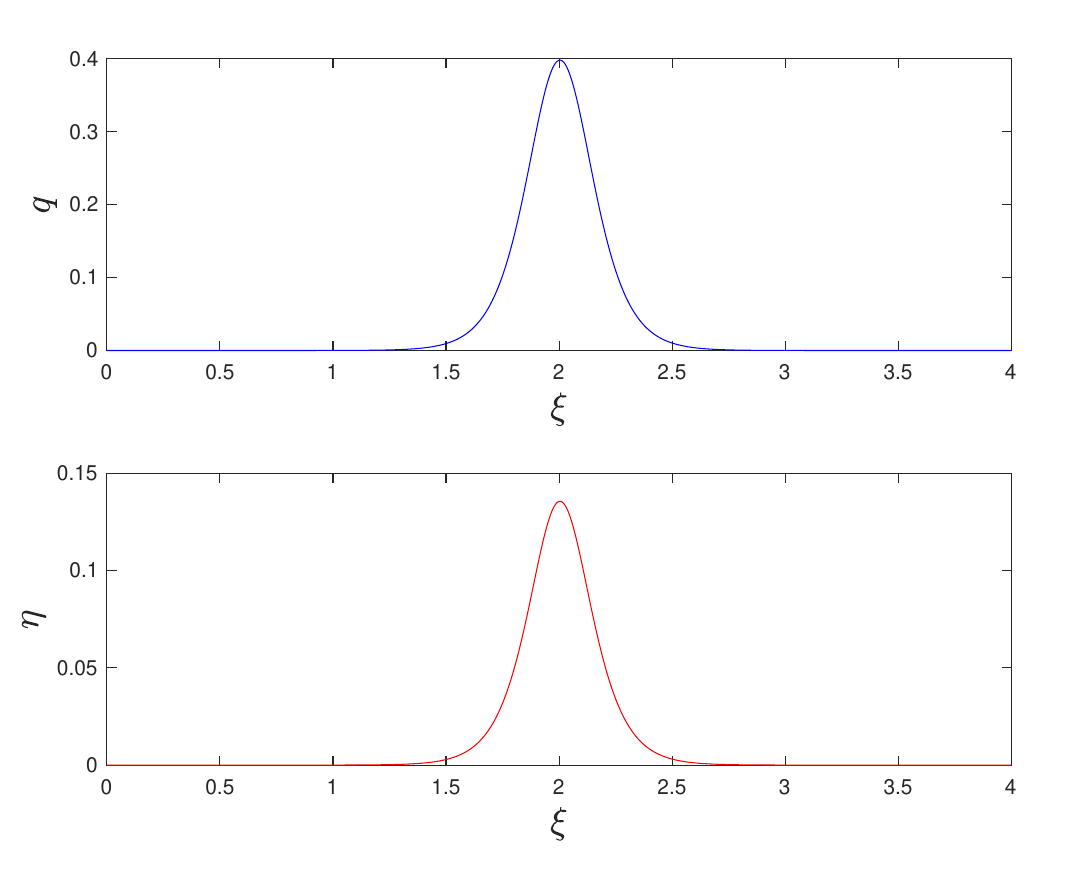}
    \caption{Construction of the traveling waves. Left panel: potential $U(q)$ corresponding to $V=10/3$ (blue continuous line). If total energy is low enough and the particle is initially in the potential well, then the orbits are periodic (black line between points $A$ and $B$). As the energy increases approaching zero from below, the period of the oscillations tends to infinity, and the trajectory becomes a traveling wave. Positive energy corresponds to open orbits. The middle panel shows the lines with constant total energy corresponding to the same potential. The thick red line is the separatrix. \revB{The right panel is obtained by numerically integrating \eqref{eq:q''} along the separatrix trajectory.}}
    \label{fig:potential}
\end{figure}

An example of potential, trajectories and traveling waves  is illustrated in Figure \ref{fig:potential}.
The first panel shows both the potential $U(q)$ corresponding to $V = 10/3$ (blue continuous line) and its best fit approximation with a cubic polynomial (magenta dashed line).  \revB{The middle panel shows the structure of this dynamical system, with two equilibria and a homoclinic connection.  The right panel is obtained by integrating \eqref{eq:q''} written as a pair of first-order equations, along the homoclinic connection, starting from a very small perturbation of the saddle equilibrium state, with the perturbation in the direction of the eigenvector of the linearized dynamical system corresponding to the positive (unstable) eigenvalue.}

\subsection{Mean profile}
\begin{figure}
     \centering
     \subfloat[Mean surface height versus $x$ for a train of solitary waves.\label{fig:three_sol}]{\includegraphics[width=\textwidth]{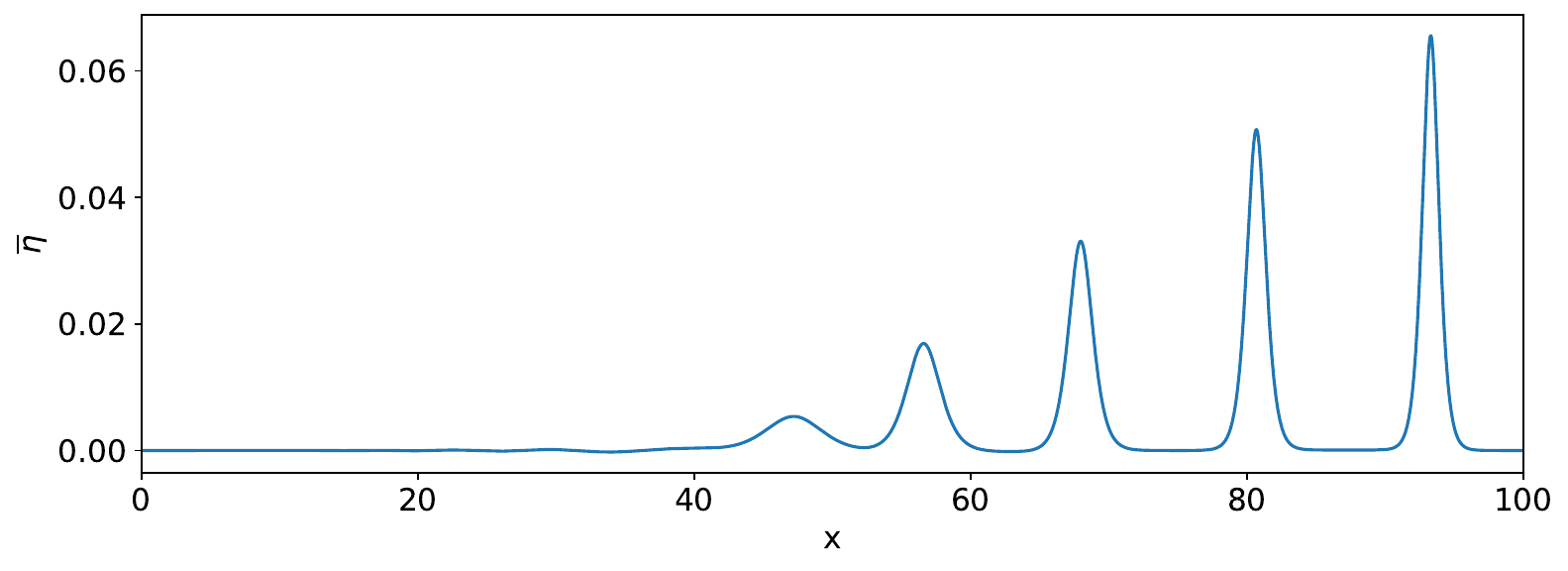}}\hfill
     \subfloat[Largest 3 waves rescaled and compared with a fitted $\sech^2$ curve.\label{fig:rescaled}]{\includegraphics[width=0.5\textwidth]{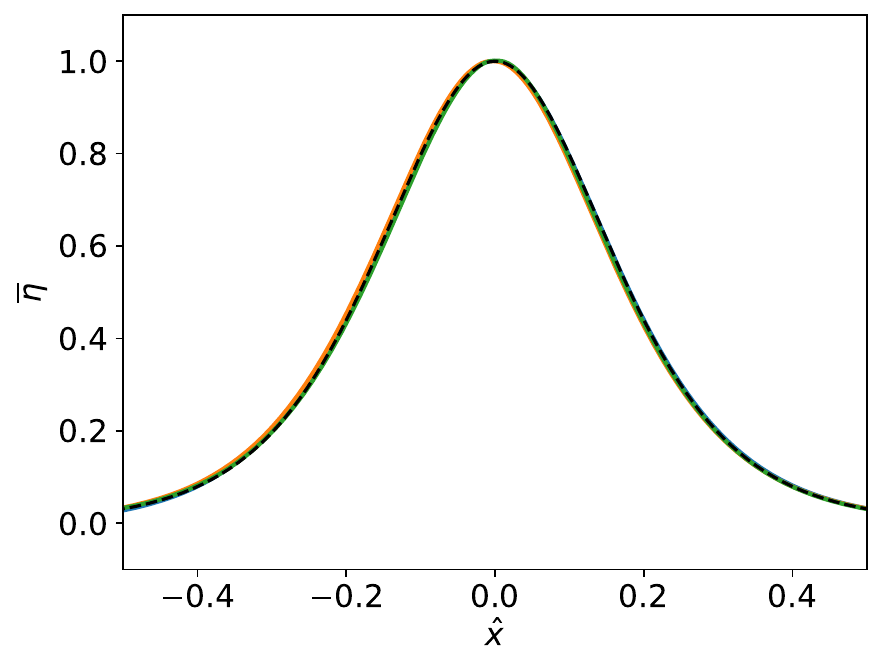}}\hfill
     \caption{The mean surface height for small-amplitude solitary waves (solid lines) is very close to $\sech^2$ (dashed line), and the waves' width scales inversely with the square root of the amplitude.}
\end{figure}

First we consider the shape of the $y$-averaged surface.  A typical $y$-averaged solution is shown in Figure \ref{fig:three_sol}.
As expected, these waves have a shape very close to the typical $\sech^2$, and seem to scale in the same way as other such solitons.  In Figure \ref{fig:rescaled} we plot each of the three tallest waves, after shifting the peak to be at $x=0$, rescaling the amplitude to 1 and rescaling the width by the square root of the amplitude.
We see that the waves very nearly coincide with the reference $\sech^2$ curve.
This is not surprising, given that the potential in the first panel of Figure \ref{fig:potential} is very well approximated by a cubic polynomial.

We have observed that much larger solitary waves have a more sharply-peaked shape; investigation of larger-amplitude solutions is the subject of ongoing work.

\subsection{Full shape}
\revA{Our traveling wave analysis shows that small-amplitude solitary wave solutions have the $y$-averaged surface profile
\begin{align}
    \mean{\eta}(x,t) \approx A \sech^2(\alpha\sqrt{A}(x-Vt),
\end{align}
where $A$ is the wave amplitude, $V$ the velocity, and $\alpha$ is related to the width; for the present scenario $\alpha \approx 4.85$.
The full approximate shape of these waves is then predicted by the formulas in Section \ref{sec:2d-structure}.
In particular, for a fixed time $t$ we have (up to translation)
\begin{align} \label{eta-full}
    \eta(x,y,t) \approx A \sech^2(\alpha\sqrt{A}x) - \delta^2 \fluctint{H^{-1}\fluctint{H}} \frac{d^2}{dx^2}\left(A\sech^2(\alpha\sqrt{A}x)\right).
\end{align}
}

\begin{figure}
     \centering
     \subfloat[Computed FV surface height (in black) versus $x$ compared to \eqref{eta-full}, for $y=-19/80$
     (dashed blue line) and $y=19/80$ (dashed orange line).  
     \label{fig:eta_y_comparison}]{\includegraphics[width=0.45\textwidth]{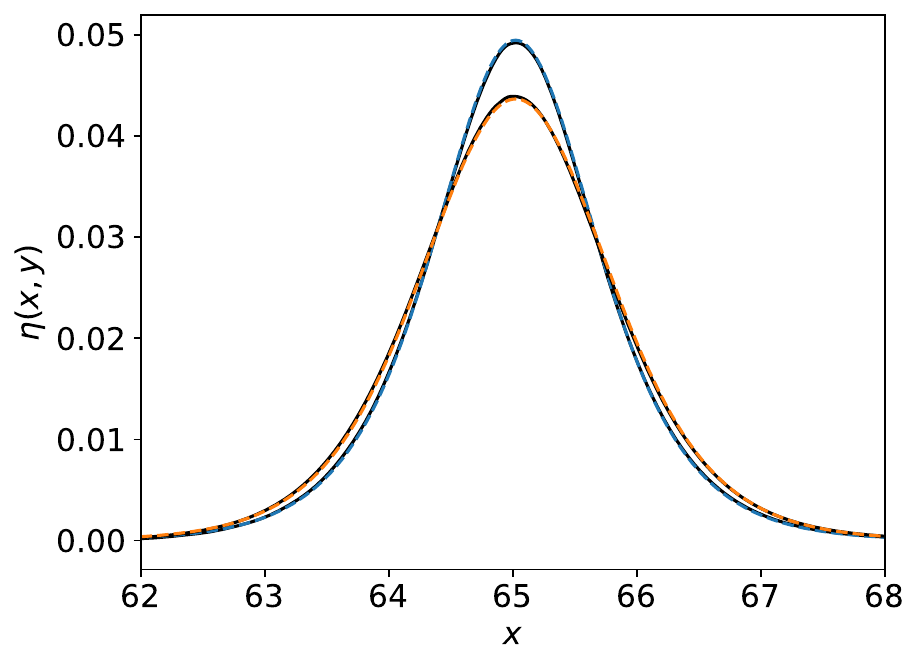}}\hfill
     \subfloat[Computed FV $y$-momentum ($p$, in black) versus $x$ compared to \eqref{p-full}, for $y=-1/80$ (dashed blue line) and $y=39/80$ (dashed orange line).\label{fig:p_y_comparison}]{\includegraphics[width=0.45\textwidth]{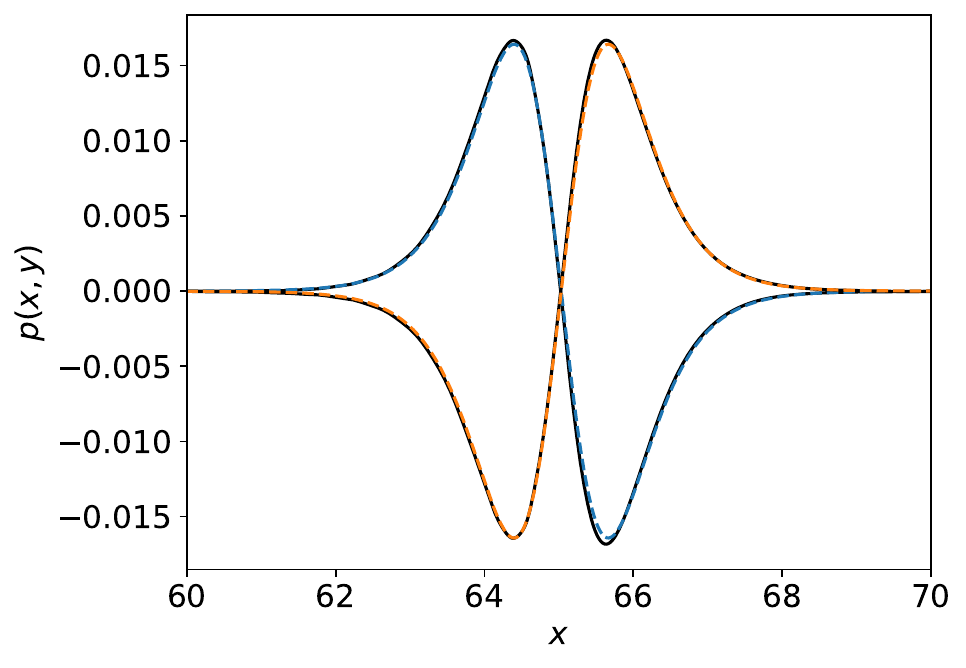}}
     \caption{Comparison of solitary wave shape \revA{(computed via FV)} with the formulas obtained from homogenization.
     \revA{FV solutions} are shown as solid black lines and predictions from homogenization are shown as colored dashed lines.
     }
     \label{fig:shape_vs_x}

\end{figure}

\begin{figure}
     \centering
     \subfloat[Computed FV surface height (in black) versus $y$ compared to \eqref{eta-full} (dashed blue line), for $x=65.1375$ .  
     \label{fig:eta_vs_y_comparison}]{\includegraphics[width=0.45\textwidth]{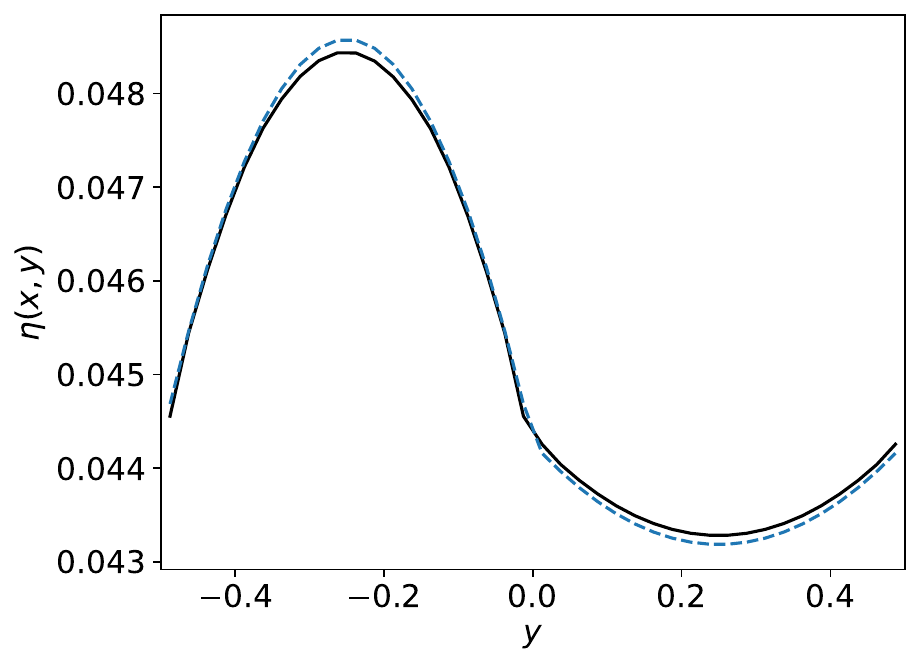}}\hfill
     \subfloat[Computed FV $y$-momentum ($p$) (in black) versus $y$ compared to \eqref{p-full}, for $x=65.4875$ (dashed blue line).\label{fig:p_vs_y_comparison}]{\includegraphics[width=0.45\textwidth]{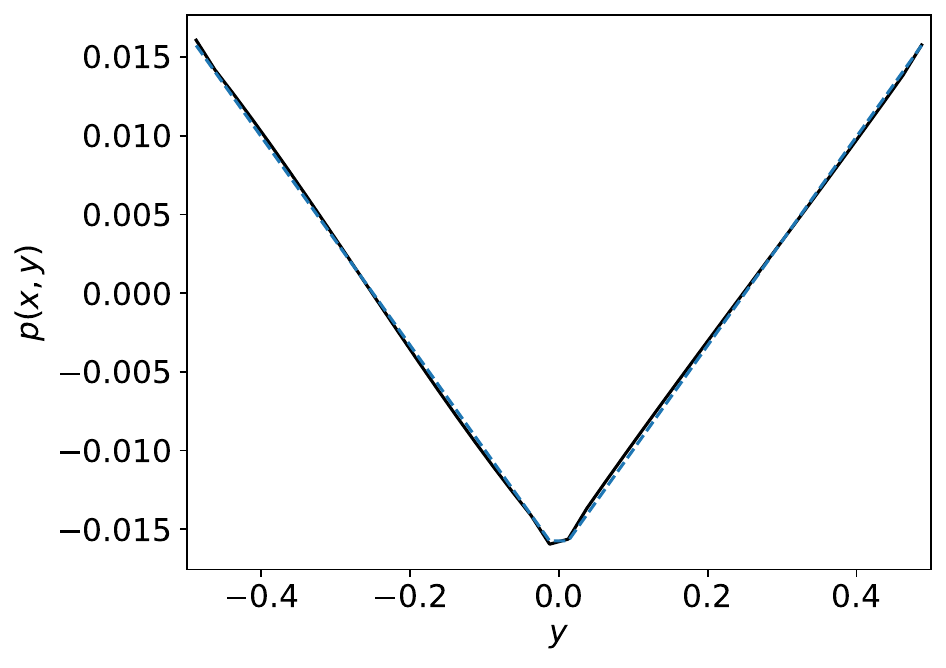}}
     \caption{\revA{Comparison of solitary wave shape (computed via FV) with the formulas obtained from homogenization, as function of $y$.
     FV solutions are shown as solid black lines and predictions from homogenization are shown as colored dashed lines.
     At the peak of the wave, $p$ vanishes, so the plot of $p$ is taken at a point away from the wave peak.}
     } \label{fig:shape_vs_y}
\end{figure}

\revA{In Figures \ref{fig:shape_vs_x} and \ref{fig:shape_vs_y}, we plot a numerical solitary wave versus the formulas from Section \ref{sec:2d-structure}, as a function of $x$ for two slices in $y$, and as a function of $y$ for slices in $x$.}
Note that here, to fit the full two-dimensional solitary wave, we have used only the mean peak amplitude as a fitting parameter.  Thus the waves we observe seem to belong to a one-parameter family.

Similar investigation of solitary waves over other
bathymetric profiles (including smooth sinusoidal bathymetry) show that the waves have, to very good approximation, the shape prescribed in Section \ref{sec:2d-structure}.

\section{Conclusion} \label{sec:conclusion}
We have studied the behavior of initially-planar shallow water waves over a bottom that varies periodically in the transverse direction.  These waves are described to good accuracy by the effective Boussinesq system \eqref{homogenized-system}, and exhibit the formation of solitary waves.  Unlike solitary wave solutions of one-dimensional hyperbolic systems with periodic coefficients \cite{leveque2003,ketcheson2023multiscale}, these are true traveling waves.  The shape of small-amplitude solitary waves is close to one that can be expressed simply in terms of elementary functions, and is predicted by the equations obtained in the process of deriving the effective Boussinesq system.   \revA{The assumption of a small amplitude wave is natural since sufficiently large-amplitude waves will exhibit shocks.  Nevertheless, in numerical experiments we have successfully produced approximate traveling wave solutions that are 2-3 times larger than those shown in this work.}

Since water waves are naturally dispersive (even over a flat bottom), it is natural to ask about the behavior of water waves over periodic bathymetry when both natural dispersion and bathymetric dispersion are accounted for.  This has been studied to some extent in \cite{chassagne2019dispersive,2021_solitary}; a full analysis starting from a dispersive water wave model
is the subject of future work, and seems to require techniques beyond what we have used herein.

Many other questions about the behavior of these waves remain open.  For instance, large-amplitude solitary waves have a different shape, and sufficiently large initial data leads to wave breaking, but the behavior of waves near the boundary between the dispersion-dominated and nonlinearity-dominated regime is complicated.  The interaction of colliding solitary waves and the behavior of periodic traveling waves in this system are also of interest.

\section*{Acknowledgment}
This work was supported by funding from King Abdullah University of Science and Technology (KAUST).  It was carried out in large part while the second author was a visiting professor at KAUST.
 G.~Russo would like to thank the Italian Ministry of University and Research (MUR) to support this research with funds coming from PRIN Project 2022 (No. 2022KA3JBA entitled  ``Advanced numerical methods for time dependent parametric partial differential equations with applications''.

\section*{Competing Interests}
The authors have no competing interests.

\bibliographystyle{plain}
\bibliography{refs}

\end{document}